\documentclass[11pt,leqno]{article}

\usepackage[T2A]{fontenc}
\usepackage{amsfonts}
\usepackage{amsmath}
\usepackage{amssymb}
\usepackage{comment}
\usepackage[dvips]{graphicx}
\usepackage{subfigure}

\newtheorem{theorem}{Theorem}[section]
\newtheorem{lemma}[theorem]{Lemma}
\newtheorem{corollary}[theorem]{Corollary}

\newcommand{\set}[1]{\left \{ #1 \right \}}                     
\newcommand{\setst}[2]{\left \{ #1 \mid #2 \right \}}           

\providecommand{\R}{\mathbb{R}}

\renewcommand{\tilde}{\widetilde}

\renewcommand{\bar}{\overline}

\newcommand{\calF}{\mathcal{F}}
\newcommand{\calD}{\mathcal{D}}

\newenvironment{proof}{\par\noindent%
{\bf Proof.\par\nopagebreak}}{\unskip\nobreak\enskip$\square$\par\bigskip}

\newenvironment{numitem}{\refstepcounter{equation}\begin{enumerate}%
\item[(\arabic{equation})]$\quad$}{\end{enumerate}}

\newcommand{\refeq}[1]{(\ref{eq:#1})}                
\newcommand{\reffig}[1]{Fig.~\ref{fig:#1}}           

\newcommand{\Xcomment}[1]{}

\begin{document}

\title
{
   Minimum Mean Cycle Problem in \\
   Bidirected and Skew-Symmetric Graphs
}

\author
{
   Maxim A. Babenko
   \thanks{
      Dept. of Mechanics and Mathematics,
      Moscow State University, Vorob'yovy Gory, 119899 Moscow,
      Russia, \textsl{email}: mab@shade.msu.ru.
      Supported by RFBR grants 03-01-00475, 05-01-02803, and 06-01-00122.
   }
   \and
   Alexander V. Karzanov
   \thanks{
      Institute for System Analysis, 9, Prospect 60 Let Oktyabrya,
      117312 Moscow, Russia, {\sl email}: sasha@cs.isa.ru. Supported by
      NWO--RFBR grant 047.011.2004.017 and RFBR grant 05-01-02805
      CNRSL\_a.
   }
}

\maketitle

 \begin{abstract}
The problem of finding, in an edge-weighted
bidirected graph $G=(V,E)$, a cycle with minimum mean weight of its
edges generalizes similar problems for both directed and undirected
graphs. (The problem is considered in two variants: for the cycles
without repeated edges and for the cycles without repeated nodes.) In
this note we develop an algorithm to solve this problem in
$O(V^2 \min(V^2, E\log V))$-time (to compare: the complexity of an improved
version of Barahona's algorithm for undirected cycles is $O(V^4)$).
Our algorithm is based on a certain general approach to minimum mean
problems and uses, as a subroutine, Gabow's algorithm for the minimum
weight 2-factor problem in a graph. The problem admits a reformulation
in terms of regular cycles in a skew-symmetric graph.
\end{abstract}

\medskip
\noindent
\emph{Keywords}: bidirected graph, skew-symmetric graph,
minimum mean cycle.

\medskip
\noindent
\emph{AMS Subject Classification}: 05C38, 05C85, 90C27.

\section{Introduction}
\label{sec:intro}

Among a variety of discrete optimization problems, an interesting
class is formed by the problems consisting in finding an object
(subset) of a given type in which the mean weight of an element is
minimized. One of the most popular problems of this sort is the
minimum mean cycle problem in an edge-weighted directed graph. The
classical algorithm due to Karp~\cite{kar-78}, based on a dynamic
programming approach, finds such a cycle in $O(nm)$ time. (Here and
later on $n$ and $m$ denote the numbers of nodes and edges,
respectively, in the input graph.) An analogous problem for cycles of
an undirected graph can also be solved efficiently:
Barahona~\cite{bar-93} reduces it to a series of
minimum weight $\emptyset$-join computations (see
also~\cite[Sec 29.11]{sch-03}), and an improved version of his algorithm
requires only $O(n)$ such computations.

This note presents an algorithm with complexity
$O(n^2 \min(n^2, m\log n))$ for a more general problem, namely, for
finding a minimum mean cycle in an edge-weighted bidirected
graph. (The concept of bidirected graph was introduced by Edmonds and
Johnson~\cite{EJ-70}; more about bidirected graphs can be found in, e.g.,
\cite{sch-03}.)

Recall that in a \emph{bidirected} graph $G=(V,E)$ edges of
three types are allowed:
a usual directed edge, or an \emph{arc}, that leaves one node and
enters another one; an edge directed \emph{from both} of its ends;
and an edge directed \emph{to both} of its ends.
When both ends of an edge coincide, the edge becomes a loop;
to slightly simplify our description we admit no loop entering and
leaving its end node simultaneously (this will lead to no loss of
generality in what follows).

A \emph{walk} in $G$ is an alternating sequence $P = (s = v_0, e_1,
v_1, \ldots, e_k, v_k = t)$ of nodes and edges such that each edge
$e_i$ connects nodes $v_{i-1}$ and $v_i$, and for $i = 1, \ldots,
k-1$, the edges $e_i,e_{i+1}$ form a \emph{transit pair} at $v_i$,
which means that one of $e_i,e_{i+1}$ enters and the other
leaves~$v_i$.
(Note that $e_1$ may enter $s$ and $e_k$ may leave $t$.)
If, in addition, $k\ge 1$, $v_0=v_k$ and the pair $e_1,e_k$ is transit
at $v_0$, $P$ is called a \emph{cycle}.
A cycle is called \emph{node-simple} if all its nodes (except
$v_0=v_k$) are different, and \emph{edge-simple} if all its edges are
different. Note that $G$ may contain edge-simple cycles but no
node-simple ones.

Our main problem is stated as follows:
  \begin{itemize}
\item[(P)] {\em Given a bidirected graph $G=(V,E)$ and a function
$w:E\to\R$ of weights of edges, Find an edge-simple cycle $C$ of
$G$ that minimizes the value $\bar w(C):=w(C)/|C|$.}
   \end{itemize}
Here $w(C)$ denotes the total weight of edges, and $|C|$ the
number of edges (the {\em length}) of $C$. So we deal with the {\em
minimum mean edge-simple cycle problem} in an edge-weighted bidirected
graph. Instead, one can state, in a similar way, the minimum mean {\em
node-simple} problem in $(G,w)$, which also may find reasonable
applications. However, the latter is easily reduced to (P) (and
therefore, can be excluded from consideration in what follows).
Indeed, modify the graph as follows: replace each node $v$ by a pair
of nodes $v_1$, $v_2$, connecting them by an edge with zero weight
from $v_1$ to $v_2$, and each edge of the original graph that enters
(leaves) $v$ make entering $v_1$ (resp. leaving $v_2$). The obtained
graph has $2n$ nodes and $m+n$ edges, all edge-simple cycles $C'$ in
it are node-simple and just correspond to the node-simple cycles $C$
of the original graph; moreover, for corresponding $C'$ and $C$, one
has $w(C')=w(C)$ and $|C'|=2|C|$.

In fact, the minimum mean cycle problem for directed graphs $G$
becomes equivalent to (P) (as well as to its node-simple version) when
$G$ is considered as a bidirected graph. As for the undirected case, an
undirected graph $G$ can be turned into a bidirected one by assigning
the orientation of each edge from both of its ends and by adding, for
each node $v$, a loop with zero weight that enters $v$ (twice). This
reduces the undirected minimum mean cycle problem to problem~(P) (but
not to its node-simple version).

\Xcomment{
\begin{figure}[tb]
   \centering
   \includegraphics{pics/reductions.1}%
   \caption{Reduction from node-simple case to edge-simple one.}
   \label{fig:node_to_edge}
\end{figure}
}

We show the following.
  \begin{theorem} \label{th:main}
Problem~(P) can be solved in $O(n^2 \min(n^2, m\log n))$ time.
  \end{theorem}

Another class of nonstandard graphs is formed by so-called
skew-symmetric graphs~--- directed graphs with involutions on the
nodes and on the arcs that change the orientation of each arc.
(This class of graphs, under the name of {\em antisymmetrical
digraphs}, was introduced by Tutte~\cite{tut-67}.) There is a close
relationship between these and bidirected graphs, and problem~(P) is,
in fact, equivalent to the problem of finding a regular cycle having
minimum mean weight in a skew-symmetric graph with symmetric weights
of arcs (precise definitions will be given later). As a consequence,
we obtain an $O(n^2 \min(n^2, m\log n))$-algorithm for the latter problem.

A natural approach to the minimum mean problem for an abstract family
$\calF$ of subsets of a set consists in reducing this problem to a
series of (auxiliary) problems of finding a member of $\calF$
with minimum total weight when the weight function is
shifted by a constant. The shifts are chosen in such a way that the
cardinalities of intermediate subsets in this process monotonically
decrease. Therefore, the process is finite, provided that the
auxiliary problem is well-solvable.
In a more general approach, $\calF$ is
extended to a family $\calD$ by adding certain sets represented as the
disjoint union of members of $\calF$. A clever choice of $\calD$ may
simplify the auxiliary problem significally, yielding an efficient
algorithm for the original problem. Just this idea is applied
in~\cite{bar-93} where undirected cycles are extended to
$\emptyset$-joins.
(One more approach to a certain class of minimum mean problems and
related topics are discussed in~\cite{kar-85}.)

In our case, the role of extended family $\calD$ plays the family of
circulations that traverse every node of a given bidirected graph at
most twice. Then the auxiliary problem becomes equivalent to the
minimum weight 2-factor problem in a certain associated undirected
graph. The latter can be solved in $O(n \min(n^2, m\log n))$ time by use of
Gabow's algorithm~\cite{gab-83}. This yields an algorithm with the
desired complexity, taking into account that the number of iterations
in the process is $O(n)$.

The above-mentioned general approach is described in
Section~\ref{sec:gen} and a proof of Theorem~\ref{th:main} is given
in Section~\ref{sec:bd}. The related problem for skew-symmetric graphs
is discussed in Section~\ref{sec:ss}.


\section{A General Approach}
\label{sec:gen}

In a general setting, the minimum mean problem is formulated as
follows:
  \begin{itemize}
\item[(M)] {\em Given a family $\calF$ of nonempty subsets of a finite
set $E$ and a function $w:E\to\R$, Find $X\in\calF$ minimizing
$\bar w(X):=w(X)/|X|$,}
   \end{itemize}
where $w(X)$ denotes $\sum_{e\in X}w(e)$. For convenience the
elements of $\calF$ are called {\em feasible sets};
so the goal is to find a feasible set with minimum mean weight.

Let us call a family $\calD$ of nonempty subsets of $E$ a
{\em disjoint extension} of $\calF$ if $\calF\subseteq\calD$
and each member $X$ of $\calD$ can be represented as the disjoint
union of feasible sets, i.e.,
  \begin{equation} \label{eq:decom}
X = X_1 \sqcup \ldots \sqcup X_k \quad
\mbox{for pairwise disjoint sets $X_1, \ldots, X_k \in \calF$}.
  \end{equation}
(To emphasize that the sets involved in the union are pairwise
disjoint, we use notation $\sqcup$ rather than $\cup$.) Let $s(\calD)$
denote the maximum cardinality (size) of a member of $\calD$.

For an arbitrary disjoint extension $\calD$ of $\calF$,
problem~(M) can be reduced to at most $s(\calD)$ problems of finding a
set $X\in\calD$ having minimum weight $w'(X)$ with respect to some
other weight function $w'$, followed by one problem of constructing a
decomposition~\refeq{decom} for the final set $X$. The
idea is rather transparent and has been encountered in special cases
(cf., e.g.,~\cite{kar-85,bar-93}).

Functions $w'$ occurring in the desired reduction are obtained by
shifting $w$ by a constant. More precisely, we say that $w'$ is the
{\em shift} of $w$ by $a\in\R$ if $w'(e) := w(e) - a$ for all $e \in
E$. We also denote $w'$ by $w^a$ and call $a$ the {\em shift number}.
The method below relies on two easy properties:
  \begin{numitem}
$\bar w^a(X) = \bar w(X) - a$ for any nonempty subset~$X \subseteq E$;
  \label{eq:shift}
  \end{numitem}
  \begin{numitem}
if $X = X_1 \sqcup \ldots \sqcup X_k$ for nonempty subsets $X_i$,
then $\bar w(X) \ge \min_i \bar w(X_i)$.
  \label{eq:mu}
  \end{numitem}

\noindent
{\bf Method.}
In the beginning put $a := \max(w(e) \colon e \in E)$. Then
iteratively construct a set $X$ and update $a$ as follows.
At each iteration, take as $X$ a member of $\calD$ with
$w^a(X)$ minimum, by solving the corresponding minimum weight problem.
Put $b:=\bar w^a(X)$. If (i) $b = 0$, then the process terminates
with constructing a decomposition~\refeq{decom} of~$X$
and outputting any of its members. Otherwise (ii) update $a := a + b$
and proceed with the next iteration.

  \begin{lemma} \label{lm:method}
The above method terminates in at most $s(\calD)+1$ iterations and
outputs an optimal solution to problem (M).
  \end{lemma}
  \begin{proof}
First of all we observe that at each iteration there exists a set
$Y\in\calD$ satisfying $w^a(Y)\le 0$. Indeed, the initial choice of
$a$ guarantees this property to hold at the first iteration.
And for each iteration, we have $w^{a'}(X)=w^a(X)-b|X|=0$, where
$a$ is the shift number at the beginning of this iteration, $X$ is the
set in $\calD$ found on it, $b=\bar w^a(X)$, and $a'$ is the next
shift number $a+b$. Then the property holds for the next iteration.

Thus, each number $b$ in the process is nonpositive. Suppose
$b=0$ happens at some iteration. Then the set $X$ found on this
iteration satisfies $w^a(X)=0$, and by the minimality of $X$, we have
$w^a(Y)\ge 0$ for all $Y\in\calD$; equivalently: $0=\bar w^a(X)\le
\bar w^a(Y)$. Now~\refeq{shift} and~\refeq{mu} imply that $\bar
w(X_i)=a$ for each member $X_i$ in a decomposition of $X$, and that
$\bar w(Y)\ge a$ for any $Y\in\calF$ (taking into account that
$\calF\subseteq\calD$). So $X_i$ is an optimal solution to (M).

Finally, consider two consecutive iterations. Let $a$ and $X$ denote
the shift numbers at the beginning of the former iteration and the set
found on it, and let $a'$ and $X'$ denote similar objects on the
latter one. Suppose that the latter iteration was not the last in the
process; then $\bar w^{a'}(X')< 0$. Also $a'=a+\bar w^a(X)$ implies
$\bar w^{a'}(X')=\bar w^a(X')-\bar w^a(X)$.
Therefore,
  $$
      w^a(X') / |X'| < w^a(X) / |X|.
  $$
This inequality is possible only if $|X'|<|X|$, in view of
$w^a(X) \le w^a(X')$ and $w^a(X)<0$.
Thus, the process is finite and the number of iterations does not
exceed $s(\calD)+1$.
  \end{proof}

\noindent
{\bf Remark.}
It is seen from the above proof that the number of iterations in the
method is estimated via the minimum cardinality of a set $X\in\calD$
minimizing $w^a(X)$ on the first iteration. Also the method can start
with any initial shift number $a$ for which one guarantees the
existence of a set $Y\in\calD$ with $w^a(Y)\le 0$ (provided that such
an $a$ can be computed efficiently).

\medskip
To illustrate the method, consider problem (M) for the family $\calF$
of simple cycles, or circuits, in a undirected graph $G=(V,E)$ with
edge weights $w$ (regarding circuits as edge sets). Since the minimum
weight problem for circuits is NP-hard when negative weights are
possible, and therefore, $\calD:=\calF$ is a bad choice, one has to
take as $\calD$ some nontrivial disjoint extension of $\calF$.
Following~\cite{bar-93}, we put $\calD$ to be the family of nonempty
$\emptyset$-joins in $G$.
(Recall that for $T\subseteq V$ with $|T|$ even, a set
$J \subseteq E$ is called a $T$-{\em join} if the set of nodes with
odd degrees in $(V,J)$ is exactly $T$.) Clearly a $\emptyset$-join
is decomposed into pairwise (edge-) disjoint circuits;
such a decomposition is constructed in
linear time. Each iteration of the method in our case consists in
finding a minimum weight $\emptyset$-join, which can be carried out in
$O(n^3)$ time (see~\cite{sch-03}).
Also the number of iterations is $O(m)$. Hence the obtained algorithm
runs in $O(n^3m)$ time.

In an improved version of this algorithm, the number of iterations
reduces to $O(n)$, yielding the overall time bound $O(n^4)$.
This relies on the existence (and the possibility to efficiently
construct) an initial shift number $a$ such that there exists a
nonempty $\emptyset$-join $J$ whose weight $w^a(J)$ is minimum and
whose cardinality is $O(n)$, and moreover, this weight is nonpositive;
cf. Remark above.

To compute the desired $a$, we order the edges of $G$ by nondecreasing
their weights:
  $$
   w(e_1) \le w(e_2) \le \ldots \le w(e_m).
  $$
We find the minimum number $k$ such that the graph
$G_k=(V,\{e_1,\ldots,e_k\})$ contains a circuit $C$. Put
$a :=\bar w(C)$; then $w^a(C) = 0$. 
Define
  \begin{eqnarray*}
    E^+ := & \setst{e \in E}{w^a(e) \ge 0} \\
    E^- := & \setst{e \in E}{w^a(e) < 0}
  \end{eqnarray*}
Let $J$ be a nonempty $\emptyset$-join with $|J|$ minimum among all
$\emptyset$-joins $J'$ minimizing $w^a(J')$. The circuit $C$ is
also a $\emptyset$-join, hence $w^a(J) \le w^a(C) = 0$. We assert that
$|J|< 2n$. Indeed, since $G_{k-1}$ is acyclic and, obviously,
contains $E^-$, we have $|J\cap E^-|\le|E^-|< k\le n$.
Also $|J\cap E^+|\le n$. For otherwise the subgraph $(V, J\cap E^+)$
would contain a circuit $C'$. Then $J':=J - C'$ is a nonempty
$\emptyset$-join satisfying $w^a(J')\le w^a(J)$ and $|J'|<|J|$,
contradicting the minimality of $J$. Thus, $|J|<2n$, and $a$ is as
required.


\section{Minimum Mean Cycles in Bidirected Graphs} \label{sec:bd}

In this section we specialize the method described in
Section~\ref{sec:gen} to find a minimum mean edge-simple cycle in a
bidirected graph $G = (V,E)$ with a weighting $w \colon E \to \R$ on
the edges. For brevity we omit the adjective ``edge-simple'' in what
follows.

Let us call a subset~$X \subseteq E$ \emph{balanced} if for each
node $v$, the numbers of edges in $X$ entering $v$ and leaving $v$ are
equal (counting twice each doubly entering or doubly leaving loop at
$v$, if any). In particular, the edge-set of any cycle is balanced. If
each node is entered by at most two edges in $X$, we say that this
balanced set is {\em small}. In a similar fashion, a cycle is called
{\em small} if it passes each node at most twice, i.e., its edge-set
is a small balanced set. Balanced sets and small cycles are related by
the following property (in fact, known in literature), which is a
bidirected analog of the fact that a circulation in a digraph is
decomposed into simple cycles.

  \begin{lemma}  \label{lm:bal-small}
Each balanced set $X$ is representable as the union of pairwise
edge-disjoint small cycles (regarding a cycle as an edge set).
  \end{lemma}
  \begin{proof}
Assuming $X\ne\emptyset$, choose a maximal edge-simple walk
$P=(v_0,e_1,v_1,\ldots,e_k,v_k)$ in the subgraph of $G$ induced by
$X$. The fact that $X$ is balanced easily implies that $v_0=v_k$ and
that $e_1,e_k$ form a transit pair at $v_0$, i.e., $P$ is a cycle.

Deleting the edges of $P$ from $X$, we again obtain a balanced set. So
it follows by induction that $X$ is representable as the union of
pairwise edge-disjoint cycles. Let the number of cycles in such a
representation be as large as possible. We assert that each cycle is
small.

Indeed, suppose some cycle $C=(v_0,e_1,v_1,\ldots,e_k,v_k)$ among
these is not small, i.e., $C$ passes some node $v$ at least three
times. Since a cycle can be considered up to reversing and cyclically
shifting, we may assume that $v=v_0=v_i=v_j$ for some $0<i<j<k$, and
that $e_1$ leaves $v_0$. If for some $p\in\{i,j\}$, the edge $e_p$
enters $v_p$, then the pair $e_1,e_p$ (as well as $e_k,e_{p+1}$) is
transit at $v$, and we can split $C$ into two cycles: the part of $C$
from $v_0$ to $v_p$ and the rest. And if $e_i$ leaves $v_i$ and $e_j$
leaves $v_j$, then the pair $e_{i+1},e_j$ is transit, and $C$ can
again be split into two cycles. A contradiction.
  \end{proof}

It is straightforward to devise a linear time algorithm that
decomposes a balanced set into small cycles.

An immediate consequence of~\refeq{mu} and Lemma~\ref{lm:bal-small} is
that for any cycle $C$, there exists a small cycle $C'$ such that
$\bar w(C')\le \bar w(C)$. Therefore, instead of all cycles, we can
take as the family $\calF$ in problem (M) the collection of (the
edge-sets of) small cycles. Also the desired disjoint extension
$\calD$ of $\calF$ can be defined to be the collection of {\em small}
balanced sets. The cardinality of any small balanced set $X$ does not
exceed $2n$ since each node is incident with at most four edges in
$X$. This implies that the number of iterations of the method applied
to these $\calF$ and $\calD$ is $O(n)$, by Lemma~\ref{lm:method}.
Each iteration consists in finding a small balanced set $X$ with
minimum weight $w'(X)$ for a current weight function $w':E\to\R$, and
it remains to explain how to solve this problem.

\begin{figure}[t]
   \centering
   \includegraphics{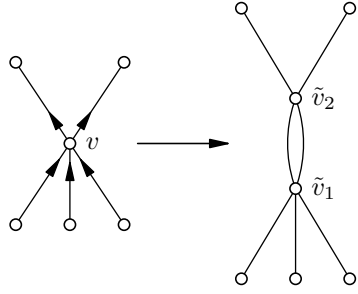}%
   \caption{Reduction to 2-factor problem.}
   \label{fig:small_balanced_set_to_2factor}
\end{figure}

We reduce it to the {\em minimum weight 2-factor problem} in a
certain undirected multigraph $\tilde G = (\tilde V, \tilde E)$,
with possible loops, formed from $G$ as follows. Each node $v \in V$
generates two nodes $\tilde v_1, \tilde v_2$ in $\tilde G$.
Each edge~$e$ in $E$ connecting nodes $u$ and $v$ generates an
edge~$\tilde e$ with the same weight connecting nodes $\tilde u_i$ and
$\tilde v_j$, by the following rule: if $e$ enters~$u$ then $i = 1$,
otherwise $i = 2$, and similarly for $v$ and $j$. In particular, a
doubly entering loop at $v$ (if any) induces an undirected loop at
$v_1$. Finally, for each $v \in V$, we add {\em two} parallel edges
connecting $\tilde v_1$ and $\tilde v_2$ in $\tilde G$ and assign them
zero weight; these edges are called \emph{auxiliary}. An example of
this transformation is depicted in~\reffig{small_balanced_set_to_2factor}.

Recall that a \emph{2-factor} in a (multi)graph is a subset of edges
such that each node of the graph is covered by exactly two edges of
this subset (counting a loop twice).
There is a natural correspondence $\phi$ between the small balanced
sets in $G$ and the 2-factors in $\tilde G$ (up to swapping parallel
auxiliary edges), and this correspondence
preserves set weights. More precisely, for a small balanced set $X$ in
$G$, $\phi(X)$ contains all edges of $\tilde G$ generated by $X$.
In addition, for each node $v \in V$ incident with $d$ entering edges
in $X$ (counting a loop twice), $\phi(X)$ contains $2-d$ auxiliary
edges connecting $\tilde v_1$ and $\tilde v_2$. One can check that
$\phi(X)$ is indeed a 2-factor. Moreover, any 2-factor of $\tilde G$
can be obtained in this way.

Applying the $O(n \min(n^2, m \log n))$-algorithm due to Gabow~\cite{gab-83} to
find a minimum weight 2-factor in $\tilde G$, we conclude that each
iteration in our method can be performed with a similar time bound.
Since the number of iterations is $O(n)$, we obtain an algorithm
to solve the original problem (P) in $O(n^2 \min(n^2, m \log n))$, as required in
Theorem~\ref{th:main}.


\section{Skew-Symmetric Graphs} \label{sec:ss}

The above result on bidirected graphs can be translated into the
language of skew-symmetric graphs. (About skew-symmetric graphs and
problems on them, see, e.g.,~\cite{tut-67,GK-96,GK-04}.)

A \emph{skew-symmetric graph} is a digraph $G=(V,E)$, with possible
multiple arcs, endowed with
two bijections $\sigma_V, \sigma_E$ such that: $\sigma_V$ is
an involution on the nodes (i.e., $\sigma_V(v)\ne v$ and
$\sigma_V(\sigma_V(v)) = v$ for each node~$v$); $\sigma_E$ is an
involution on the arcs; and for each arc $a$ from $u$ to $v$,
$\sigma_E(a)$ is an arc from $\sigma_V(v)$ to $\sigma_V(u)$. For
brevity, the mappings $\sigma_V, \sigma_E$ are combined into one
mapping $\sigma$ on $V\cup E$, which is called the
\emph{symmetry} (rather than skew-symmetry) of $G$. For a node (arc)
$x$, its symmetric node (arc) $\sigma(x)$ is also called the
\emph{mate} of $x$, and we will use notation with primes for
mates, denoting $\sigma(x)$ by $x'$.

Observe that if $G$ contains an arc $a$ from a node $v$
to its mate $v'$, then $a'$ is also an arc from $v$ to $v'$.

The symmetry $\sigma$ is extended in a natural way to walks, cycles
and other objects in $G$. In particular, two walks or
cycles are symmetric to each other if the elements of one of them
are symmetric to those of the other and go in the reverse order:
for a walk $P = (v_0, a_1, v_1, \ldots, a_k, v_k)$, the symmetric
walk $\sigma(P)$ is $(v'_k, a'_k, v'_{k-1}, \ldots, a'_1, v'_0)$.

Next we explain a relationship between skew-symmetric and
bidirected graphs (cf.~\cite[Sec.~2]{GK-04}).
Given a skew-symmetric graph $G = (V,E)$, choose an arbitrary
partition $\pi=\set{V_1, V_2}$ of $V$ such that $V_2$ is symmetric to
$V_1$. Then $G$ and $\pi$ determine the bidirected graph $\bar G =
(\bar V, \bar E)$ with $\bar V := V_1$ whose edges correspond to the
pairs of symmetric arcs in $G$. More precisely, arc mates $a,a'$ of
$G$ generate one edge $e$ of $\bar G$ connecting nodes $u, v \in V_1$
such that: (i) $e$ goes from $u$ to $v$ if one of $a, a'$ goes from
$u$ to $v$ (and the other goes from $v'$ to $u'$ in $V_2$); (ii)
$e$~leaves both $u,v$ if one of $a,a'$ goes from $u$ to $v'$ (and the
other from $v$ to $u'$); (iii) $e$~enters both $u, v$ if one of $a,a'$
goes from $u'$ to $v$ (and the other from $v'$ to $u$). In particular,
$e$ is a loop if $a, a'$ connect a pair of symmetric nodes.

\begin{figure}[tb]
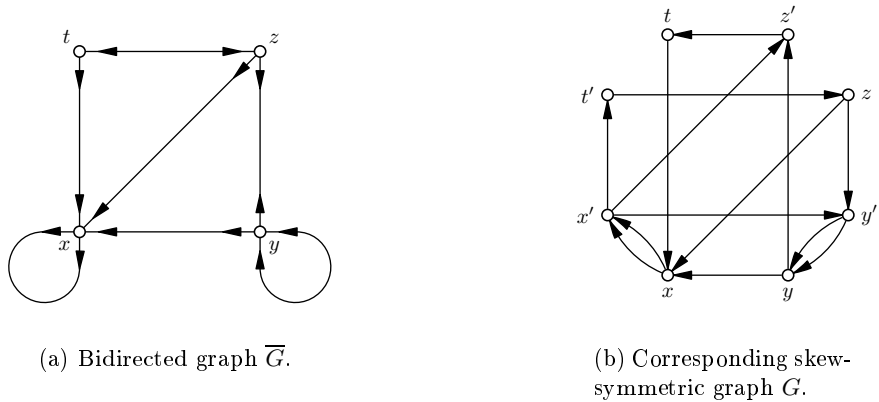

    \centering
    \subfigure[Bidirected graph~$\bar G$.]{
      \includegraphics{pics/examples.1}%
    }
    \hspace{3cm}%
    \subfigure[Corresponding skew-sym\-met\-ric graph~$G$.]{
      \includegraphics{pics/examples.2}%
    }
    \caption{Related bidirected and skew-symmetric graphs.}
    \label{fig:sk-bi}
\end{figure}

Conversely, a bidirected graph $\bar G = (\bar V, \bar E)$
determines a skew-symmetric graph $G = (V,E)$ with symmetry $\sigma$
as follows. Take a copy $\sigma(v)$ of each element $v$ of $\bar V$,
forming the set $\bar V' := \setst{\sigma(v)}{v\in \bar V}$. Now set
$V := \bar V \sqcup \bar V'$. For each edge $e$ of $\bar G$ connecting
nodes $u$ and $v$, assign two ``symmetric'' arcs $a,a'$ in $G$ so as
to satisfy (i)--(iii) above (where $u' = \sigma(u)$ and $v' =
\sigma(v)$). An example is depicted in~\reffig{sk-bi}.

(Note that one bidirected graph generates one skew-symmetric graph, by
the second construction. On the other hand, one skew-symmetric graph
generates a set of bidirected graphs, depending on the partition $\pi$
of $V$, by the first construction. One can check that among these
bidirected graphs, one is obtained from another by the following
operation: choose a subset $X$ of nodes, and for each node $v$ in
$X$ and each edge $e$ incident with $v$, reverse the direction of $e$
at $v$. This implies that the sets of walks (cycles) in these
bidirected graphs are the same.)

There is a one-to-one correspondence between walks in $\bar G$ and
(directed) walks in $G$. More precisely, let $\tau$ be the natural
mapping of $V \cup E$ to $\bar V \cup \bar E$ (obtained by identifying
the pairs of symmetric nodes and arcs). Each walk $P = (v_0, a_1, v_1,
\ldots, a_k, v_k)$ in $G$ (where $a_i$ is an arc from $v_{i-1}$ to
$v_i$) induces the sequence
  $$
   \tau(P) := (\tau(v_0), \tau(a_1), \tau(v_1), \ldots,
 \tau(a_k), \tau(v_k))
$$
of nodes and edges in $\bar G$. One can see that $\tau(P)$ is a walk
in $\bar G$ and that $\tau(P')$ is the walk reverse to $\tau(P)$.
Moreover, for any walk $\bar P$ in $\bar G$, there is exactly one
preimage $\tau^{-1}(\bar P)$.

Finally, following terminology of~\cite{GK-96}, a walk (cycle) in a
skew-symmetric graph is called \emph{regular} if it is arc-simple and
contains no pair of symmetric arcs (while pairs of symmetric nodes in
it are allowed). Certain problems on regular walks and cycles are
studied in~\cite{tut-67,GK-96} (and some other works). One more
problem is:
  \begin{itemize}
\item[(S)] {\em Given a skew-symmetric graph $G=(V,E)$ and a symmetric
weight function $w:E\to\R$ (i.e., $w(a)=w(a')$ for all $a\in E$),
find a regular cycle $C$ in $G$ whose mean weight $\bar w(C)$ is
minimum. }
  \end{itemize}

It is not difficult to check that for each regular cycle $C$ in $G$,
the corresponding cycle $\tau(C)$ in the bidirected graph $\bar G$ is
edge-simple, and vice versa. Moreover, the weights (and the mean
weights) of the corresponding cycles are equal, assuming that for
$e\in E$, the weight of $\tau(e)$ is defined to be $w(e)$ (this is
well-defined since $w$ is symmetric). Therefore, problem (S) is
equivalent to (P), and we obtain the following.

  \begin{corollary}
Problem (S) can be solved in $O(n^2 \min(n^2, m\log n))$ time.
  \end{corollary}

\bibliographystyle{amsplain}
\bibliography{main}

\end{document}